\documentclass[12pt]{amsart}
\usepackage{amsmath,amssymb}
\usepackage{amscd}
\usepackage{mathrsfs}
\usepackage[all]{xy}
\usepackage{pstricks}
\UseComputerModernTips
\usepackage{graphicx}
\DeclareGraphicsExtensions{.eps} {\bf }

\theoremstyle{plain}

\newtheorem{lemma}{Lemma}[section]

\newcommand{\edge}[1]{\ar@{-}[#1]}
\theoremstyle{definition}
\newtheorem{definition}{Definition}[section]
\newtheorem{example}{Example}[section]
\theoremstyle{remark}

\def\rank{\operatorname{rank}}

\def\span{\operatorname{span}}

\def\P{\mathcal     P}

\numberwithin{equation}{section} \numberwithin{equation}{section}
\textheight=8.75in \theoremstyle{plain}
\begin{document}
\title[On the msr  of a given graph of at most seven vertices]{Finding msr of a given graph of at most seven vertices by giving vector representations}
\author{Xinyun Zhu}
\address{Department of Mathematics\\  University of Texas of Permian Basin\\ Odessa, TX 79762}
\email{zhu\_x@utpb.edu} \date{\today}
\begin{abstract} 

In this paper, we study the  minimum rank among positive semidefinite matrices with a given graph of at most seven vertices (msr) by giving vector representations.
\end{abstract}
\maketitle
{\\

{} 
\section{Introduction}
Given a connected graph $G$ with $n$ vertices, we associate to $G$ a set $\mathcal{H}(G)$  of Hermitian $n\times n$ matrices by the following way,
\[\mathcal{H}(G)=\{ A|A=A^{\star}, a_{ij}\ne 0  \quad \text{for}\quad  i\ne j  \quad \text{if and only if } (i,j) \quad \text{ is an edge of}\quad G\}\]
where $A^{\star}$ is the complex conjugate of $A$ and $a_{ij}$ is the $ij$-entry of $A.$
 We define $\mathcal{P}(G)$ be the subset of $\mathcal{H}(G)$ whose members are  positive semidefinite matrices. For a given graph $G$, $\mathcal{P}(G)$ is non-empty because laplacian matrix $L(G)\in \P(G).$
Define
\[msr(G)=\min_{A\in \mathcal{P}(G)}\rank(A)\]

To introduce the result of this paper, we need the following preparation.

\begin{definition}\label{def1}An induced subgraph $H$ of a graph $G$ is obtained by
deleting all vertices except for the vertices in a subset $S$. For a graph $G$, we consider its "tree size," denoted $ts(G)$, which is the number of
vertices in a maximum induced tree.
\end{definition}
\begin{lemma}\cite[Lemma 2.1]{BH}\label{lemma21} If $H$ is an induced subgraph of a connected graph $G,$ then $msr(H)\le msr(G).$
\end{lemma}
Given a simple connected graph $G,$ we can define the tree size $ts(G)$ of $G.$  It is known that for a tree $T,$ $msr(T)=ts(T)-1.$ Hence by Lemma~\ref{lemma21},we get $ts(G)-1\le msr(G).$ 
\begin{lemma}\cite[Corollary 3.5]{BH} If a simple connected graph $G$ has a pendant vertex $v,$ which is simply vertex of degree $1,$  then $msr(G)=msr(G-v)+1.$
\end{lemma}

\section{Finding msr by giving vector representations}\label{examples}

Suppose $G$ is a connected graph with vertex set $V(G) = \{v_1, v_2, \cdots , v_n\}.$ We call
a set of vectors
$V = \{\vec{v_1},\cdots, \vec{v_n}\}$ in $\mathbb{C}_m$
a vector representation (or orthogonal
representation) of $G$ if  for $i\not= j\in V(G),$ $\langle\vec{v_i},\vec{v_j}\rangle\not=0$ whenever  $i$  and $j$ are adjacent in $G$ and $\langle\vec{v_i},\vec{v_j}\rangle=0$ whenever  $i$  and $j$ are not adjacent in $G$. For any matrix $B,$ $A=B^{\star}B$ is a semi-definite positive matrix and $\rank(A)=\rank(B)$. For any semi-definite positive matrix $A,$ there exists a matrix $B$ such that $A=B^{\star}B$ and $\rank(A)=\rank(B).$

For a large family of graphs, \cite{BH} has given  a good  way to get msr.  But \cite{BH} did not give the msr problem completely.  In this section we solve the msr for graphs with at most seven vertices but does not satisfies the conditions of graphs in \cite{BH}.   For the readers convenience, we list all these graphs here: $G706$, $G710$, $G817$, $G836$, $G864,$ $G867$, $G870,$ $G872,$ $G876$, $G877$, $G946$, $G954$, $G955$, $G979$, $G982$, $G992$, $G997-G1000$, $G1003-G1007$, $G1053$, $G1056$, $G1060$, $G1065$, $G1069$, $G1084$, $G1089-G1097,$ $G1100,$ $G1101, $ $1104,$ $G1105,$ $G1123,$ $G1125,$ $G1135,$ $G1145,$  $G1146$, $G1148,$ $G1149,$ $G1152-G1157,$ $G1159,$ $G1160,$ $G1165, $ $G1167,$ $G1168,$ $G1170, $ $G1176,$ $G1179,$  $G1189, $ $G1191,$ $G1194-G1197$ $G1199,$ $G1200,$ $G1202, $ $G1205,$ $G1207-G1212,$ $G1222,$  $G1224,$ $G1228$, $G1230$, $G1231$, $G1233$, $G1241$, $G1242$, $G1248,$ $G1250.$
We give the msr of all those graphs by finding the vector representations.

In the following examples,  for vectors $v_i$ and $v_j,$ by abusing notations, we write $\langle v_i, v_j\rangle=v_i\cdot v_j.$ If $v_i\cdot v_j=0,$ then we write $v_i\perp v_j.$
\begin{example}
\[msr(G706)=5\]
The tree size of $G706$ is $5.$ First let's show $msr(G706)>4.$
Let $\{v_1, v_2, v_3, v_4, v_5, v_6, v_7\}$ be a vector representation of $G706.$ 
It follows that $v_1\perp \span\{v_3,v_4,v_5,v_6\},$ 
$v_2\perp \span\{v_5,v_6\},$ 
$v_3\perp \span\{v_4,v_5,v_6\},$
 and $v_7\perp \span\{v_4,v_5\}.$
Hence $v_5\notin \span \{v_6\}.$
Otherwise $v_6\cdot v_7=0.$
Similarly, $v_4\notin \span \{v_5,v_6\}.$
Otherwise $v_2\cdot v_4=0.$
Similarly,$v_3\notin \span \{v_4,v_5,v_6\}.$
Otherwise $v_3=0.$
If the $msr(G706)$ is $4$ , we get $v_1\in \span\{v_3,v_4,v_5,v_6\}.$  It follows that $v_1=0.$  This is a contradiction.

On the other hand, let

\[A=\left[ \begin{matrix} 1 &1 &0&0&0&0&1\\
1&4&1&1&0&0&2\\
0&1&1&0&0&0&1\\
0&1&0&1&1&1&0\\
0&0&0&1&2&2&0\\
0&0&0&1&2&3&1\\
1&2&1&0&0&1&3\\
\end{matrix}
\right]\]

Let \[B=\left[\begin{matrix} 1 &1 &0&0&0&0&1\\
0&1&1&0&0&0&1\\
0&1&0&1&1&1&0\\
0&-1&0&0&1&1&0\\
0&0&0&0&0&1&1\\
\end{matrix}
\right]\]
Then $B^{T}B=A$, $A\in \P(G706)$ and $\rank(B)=5.$
\end{example}

\begin{example}
\[msr(710)=5\]
The tree size of $G710$ is $5.$ First let's show $msr(G710)>4.$

Suppose $msr(G710)=4.$  Let $\{v_1,v_2, v_3, v_4, v_5, v_6, v_7\}$ be a vector representation. Then $v_1\perp \span\{v_3, v_4, v_5, v_6\},$ $v_2\perp\span\{v_4, v_5, v_6\},$ $v_3\perp\span\{v_5,v_7\},$ $v_7\perp\span \{v_4, v_5\}.$  It follows that $v_3, v_4, v_5, v_6$ must be linearly dependent.  If $v_4, v_5, v_6$ are linearly independent, then $v_3\in \span \{v_4, v_5, v_6\}.$  Hence $v_2\cdot v_3=0.$ This is a contradiction.  Hence $v_4, v_5, v_6$ are linearly dependent.  We know $v_4$ and $v_5$ are linearly independent because $v_3\cdot v_4\not=0$ and $v_3\cdot v_5=0.$ It follows that $v_6\in \span \{v_4, v_5\}.$ It follows that $v_6\cdot v_7=0.$ This is a contradiction.

On the other hand, let

\[A=\left[ \begin{matrix} 1 &1 &0&0&0&0&1\\
1&2&1&1&0&0&1\\
0&1&3&1&0&-1&0\\
0&0&1&1&1&1&0\\
0&0&0&1&2&3&0\\
0&0&-1&1&3&6&1\\
1&1&0&0&0&1&2\\
\end{matrix}
\right]\]

Let \[B=\left[\begin{matrix} 1 &1 &0&0&0&0&1\\
0&0&1&1&1&1&0\\
0&0&-1&0&1&2&0\\
0&0&0&0&0&1&1\\
0&1&1&0&0&0&0\\
\end{matrix}
\right]\]
Then $B^{T}B=A$, $A\in \P(G710)$ and $\rank(B)=5.$
\end{example}

\begin{example}
\[msr(G817)=4\]
The tree size of $G817$ is  $5$ and let
\[A=\left[ \begin{matrix} 1 &1 &0&0&0&0&1\\
1&5&2&0&0&0&3\\
0&2&2&2&0&1&0\\
0&0&2&5&1&2&-1\\
0&0&0&1&1&0&1\\
0&0&1&2&0&1&-1\\
0&3&0&-1&1&-1&4\\
\end{matrix}
\right]\]

Let \[B=\left[\begin{matrix} 1 &1 &0&0&0&0&1\\
0&2&1&0&0&0&1\\
0&0&0&1&1&0&1\\
0&0&1&2&0&1&-1\\
\end{matrix}
\right]\]
Then $B^{T}B=A$,$A\in \P(G817)$ and $\rank(B)=4.$

\end{example}

\begin{example}
\[msr(G836)=4\]
The tree size of $G836$ is $5$  and let
\[A=\left[ \begin{matrix} 1 &1 &0&0&0&0&1\\
1&2&1&0&0&0&2\\
0&1&3&1&1&3&0\\
0&0&1&1&0&1&0\\
0&0&1&0&1&2&-1\\
0&0&3&1&2&5&-2\\
1&2&0&0&-1&-2&3\\
\end{matrix}
\right]\]

Let \[B=\left[\begin{matrix} 1 &1 &0&0&0&0&1\\
0&1&1&0&0&0&1\\
0&0&1&1&0&1&0\\
0&0&1&0&1&2&-1\\
\end{matrix}
\right]\]
Then $B^{T}B=A$, $A\in \P(G836)$ and $\rank(B)=4.$

\end{example}

\begin{example}
\[msr(G864)=4\]
The tree size of  $G864$ is $5$  and let
\[A=\left[ \begin{matrix} 1 &1 &0&0&0&1&-2\\
1&5&2&0&0&-1&0\\
0&2&3&1&1&0&3\\
0&0&1&2&0&0&1\\
0&0&1&0&2&2&3\\
1&-1&0&0&2&4&0\\
-2&0&3&1&3&0&10\\
\end{matrix}
\right]\]

Let \[B=\left[\begin{matrix} 1 &1 &0&0&0&1&-2\\
0&0&1&1&1&1&2\\
0&2&1&0&0&-1&1\\
0&0&0&1&-1&-1&-1\\
\end{matrix}
\right]\]
Then $B^{T}B=A$, $A\in \P(G864)$ and $\rank(B)=4.$

\end{example}

\begin{example}
\[msr(G867)=4\]
The tree size of  $G867$ is $5$  and let

\[A=\left[ \begin{matrix} 1 &1 &0&1&0&1&0\\
1&7&5&0&-1&0&3\\
0&5&5&0&0&0&3\\
1&0&0&7&1&2&-1\\
0&-1&0&1&1&1&0\\
1&0&0&2&1&2&0\\
0&3&3&-1&0&0&2\\
\end{matrix}
\right]\]

Let \[B=\left[\begin{matrix} 1 &1 &0&1&0&1&0\\
0&1&1&-2&0&0&1\\
0&-1&0&1&1&1&0\\
0&2&2&1&0&0&1\\
\end{matrix}
\right]\]
Then $B^{T}B=A$, $A\in \P(G867)$  and $\rank(B)=4.$

\end{example}
\begin{example}
\[msr(G870)=4\]
The tree size of  $G870$ is $5$  and let

\[A=\left[ \begin{matrix} 1 &1 &0&0&0&2&-2\\
1&4&1&0&0&-2&0\\
0&1&1&1&0&0&4\\
0&0&1&2&1&1&5\\
0&0&0&1&2&-2&0\\
2&-2&0&1&-2&14&0\\
-2&0&4&5&0&0&22\\
\end{matrix}
\right]\]

Let \[B=\left[\begin{matrix} 1 &1 &0&0&0&2&-2\\
0&1&1&1&0&0&4\\
0&-1&0&1&1&1&1\\
0&-1&0&0&-1&3&1\\
\end{matrix}
\right]\]
Then $B^{T}B=A$,$A\in \P(G870)$ and $\rank(B)=4.$

\end{example}
\begin{example}
\[msr(G872)=4\]
The tree size of  $G872$ is $5$  and let

\[A=\left[ \begin{matrix} 1 &1 &0&0&1&1&0\\
1&4&1&-1&0&0&0\\
0&1&1&0&1&-1&0\\
0&-1&0&1&1&1&1\\
1&0&1&1&4&0&0\\
1&0&-1&1&0&4&2\\
0&0&0&1&0&2&2\\
\end{matrix}
\right]\]

Let \[B=\left[\begin{matrix} 1 &1 &0&0&1&1&0\\
0&1&1&0&1&-1&0\\
0&-1&0&1&1&1&1\\
0&1&0&0&-1&1&1\\
\end{matrix}
\right]\]
Then $B^{T}B=A$,$A\in \P(G872)$ and $\rank(B)=4.$

\end{example}

\begin{example}
\[msr(G876)=4\]
The tree size of  $G876$ is $5$  and let

\[A=\left[ \begin{matrix} 1 &1 &1&1&0&0&0\\
1&3&0&0&3&2&0\\
1&0&3&0&0&-1&2\\
1&0&0&6&0&-1&-3\\
0&3&0&0&9&3&1\\
0&2&-1&-1&3&2&0\\
0&0&2&-3&1&0&3\\
\end{matrix}
\right]\]

Let \[B=\left[\begin{matrix} 1 &1 &1&1&0&0&0\\
0&1&-1&1&2&1&-1\\
0&0&1&0&2&0&1\\
0&1&0&-2&1&1&1\\
\end{matrix}
\right]\]
Then $B^{T}B=A$, $A\in \P(G876)$ and $\rank(B)=4.$

\end{example}

\begin{example}
\[msr(G877)=4\]
The tree size of  $G877$ is $5$.

Let

\[A=\left[ \begin{matrix} 1 &1 &0&0&0&-1&1\\
1&3&1&0&3&0&0\\
0&1&2&1&1&0&-5\\
0&0&1&1&0&-1&-3\\
0&3&1&0&5&1&0\\
-1&0&0&-1&1&3&0\\
1&0&-5&-3&0&0&14\\
\end{matrix}
\right]\]

Let \[B=\left[\begin{matrix} 1 &1 &0&0&0&-1&1\\
0&1&1&0&1&1&-2\\
0&0&1&1&0&-1&-3\\
0&1&0&0&2&0&1\\
\end{matrix}
\right]\]
Then $B^{T}B=A$, $A\in \P(G877)$ and $\rank(B)=4.$
\end{example}
\begin{example}
\[msr(G946)=4\]
The tree size of $G946$ is $3.$  First let's show that $msr(G946)>3.$

Let $\{v_1,v_2,v_3,v_4,v_5,v_6,v_7\}$ be a vector representation of $G946.$ We have $v_1\perp\{v_3,v_4,v_5,v_7\},$ $v_2\perp \{v_4,v_5\},$ $v_3\perp\{v_6\},$ $v_4\perp\{v_6,v_7\}.$  If $msr(G946)=3,$ then $v_3=av_4+bv_7$ for some nonzero $a$ and $b.$  Hence $0=av_4\cdot v_6+bv_7\cdot v_6=bv_7\cdot v_6.$ This is a contradiction.

On the other hand, let
\[A=\left[ \begin{matrix} 1 &1 &0&0&0&-1&0\\
1&2&1&0&0&-2&1\\
0&1&3&1&2&0&3\\
0&0&1&1&1&0&0\\
0&0&2&1&2&1&2\\
-1&-2&0&0&1&3&1\\
0&1&3&0&2&1&5\\
\end{matrix}
\right]\]

Let \[B=\left[\begin{matrix} 1 &1 &0&0&0&-1&0\\
0&0&1&1&1&0&0\\
0&1&1&0&0&-1&1\\
0&0&1&0&1&1&2\\
\end{matrix}
\right]\]
Then $B^{T}B=A$, $A\in \P(G946)$ and $\rank(B)=4.$
\end{example}

\begin{example}
\[msr(G954)=4\]

The tree size of $G954$ is $4$.

Let

\[A=\left[ \begin{matrix} 1 &1 &0&0&0&1&1\\
1&15&5&1&0&-1&6\\
0&5&2&0&0&-1&2\\
0&1&0&1&1&1&0\\
0&0&0&1&3&2&0\\
1&-1&-1&1&2&3&0\\
1&6&2&0&0&0&3\\
\end{matrix}
\right]\]

Let \[B=\left[\begin{matrix} 1 &1 &0&0&0&1&1\\
0&3&1&0&-1&-1&1\\
0&1&0&1&1&1&0\\
0&2&1&0&1&0&1\\
\end{matrix}
\right]\]
Then $B^{T}B=A$, $A\in \P(G954)$ and $\rank(B)=4.$
\end{example}

\begin{example}
\[msr(G955)=4\]
The tree size of $G955$ is $5$.

Let

\[A=\left[ \begin{matrix} 1 &1 &0&0&0&1&0\\
1&6&3&0&-1&3&-2\\
0&3&3&1&0&2&0\\
0&0&1&1&0&1&1\\
0&-1&0&0&2&-1&1\\
1&3&2&1&-1&3&0\\
0&-2&0&1&1&0&2\\
\end{matrix}
\right]\]

Let \[B=\left[\begin{matrix} 1 &1 &0&0&0&1&0\\
0&0&1&1&0&1&1\\
0&1&1&0&1&0&0\\
0&2&1&0&-1&1&-1\\
\end{matrix}
\right]\]
Then $B^{T}B=A$, $A\in \P(G955)$ and $\rank(B)=4.$
\end{example}
\begin{example}
\[msr(G979)=4\]
The tree size of $G979$ is $5$.

Let

\[A=\left[ \begin{matrix} 1 &1 &0&0&0&0&1\\
1&3&1&0&0&1&2\\
0&1&1&1&0&1&1\\
0&0&1&3&-1&2&0\\
0&0&0&-1&1&-1&1\\
0&1&1&2&-1&2&0\\
1&2&1&0&1&0&3\\
\end{matrix}
\right]\]

Let \[B=\left[\begin{matrix} 1 &1 &0&0&0&0&1\\
0&1&1&1&0&1&1\\
0&0&0&-1&1&-1&1\\
0&1&0&-1&0&0&0\\
\end{matrix}
\right]\]
Then $B^{T}B=A$, $A\in \P(G979)$ and $\rank(B)=4.$
\end{example}

\begin{example}
\[msr(G982)=4\]
The tree size of $G982$ is $5$.

Let

\[A=\left[ \begin{matrix} 1 &1 &0&0&0&0&1\\
1&4&2&0&0&6&2\\
0&2&2&4&0&4&0\\
0&0&4&26&-2&2&-6\\
0&0&0&-2&2&-2&2\\
0&6&4&2&-2&14&0\\
1&2&0&-6&2&0&4\\
\end{matrix}
\right]\]

Let \[B=\left[\begin{matrix} 1 &1 &0&0&0&0&1\\
0&1&1&1&1&1&1\\
0&1&1&3&-1&3&-1\\
0&1&0&-4&0&2&1\\
\end{matrix}
\right]\]
Then $B^{T}B=A$, $A\in \P(G982)$ and $\rank(B)=4.$
\end{example}

\begin{example}
\[msr(G992)=4\]
The tree size of $G992$ is $5$.

Let

\[A=\left[ \begin{matrix} 1 &1 &0&0&0&1&2\\
1&2&1&0&0&3&1\\
0&1&2&1&0&4&0\\
0&0&1&2&1&4&0\\
0&0&0&1&1&2&-1\\
1&3&4&4&2&13&0\\
2&1&0&0&-1&0&7\\
\end{matrix}
\right]\]

Let \[B=\left[\begin{matrix} 1 &1 &0&0&0&1&2\\
0&0&1&1&0&2&1\\
0&0&0&1&1&2&-1\\
0&1&1&0&0&2&-1\\
\end{matrix}
\right]\]
Then $B^{T}B=A$, $A\in \P(G992)$ and $\rank(B)=4.$
\end{example}

\begin{example}
\[msr(G997)=4\]
The tree size of $G997$ is $5$.

Let

\[A=\left[ \begin{matrix} 1 &2 &0&0&0&1&1\\
2&5&1&0&0&3&1\\
0&1&2&1&0&2&0\\
0&0&1&2&1&0&2\\
0&0&0&1&1&-1&1\\
1&3&2&0&-1&4&0\\
1&1&0&2&1&0&4\\
\end{matrix}
\right]\]

Let \[B=\left[\begin{matrix} 1 &2 &0&0&0&1&1\\
0&0&1&1&0&1&1\\
0&0&0&1&1&-1&1\\
0&1&1&0&0&1&-1\\
\end{matrix}
\right]\]
Then $B^{T}B=A$, $A\in \P(G997)$ and $\rank(B)=4.$
\end{example}
\begin{example}
\[msr(G998)=4\]
The tree size of $G998$ is $3.$  First let's show that $msr(G998)>3.$

Let $\{v_1,v_2,v_3,v_4,v_5,v_6,v_7\}$ be a vector representation of $G998.$ We have $v_1\perp\{v_3,v_4,v_6\},$ $v_2\perp \{v_4,v_6\},$ $v_3\perp\{v_5\},$ $v_4\perp\{v_5\},$  $v_5\perp\{v_7\},$$v_6\perp\{v_7\}.$If $msr(G998)=3,$ then $av_3+bv_4+cv_6=0$ for some nonzero $a,\,b$ and $c.$  Hence $0=av_3\cdot v_5+bv_4\cdot v_5+cv_6\cdot v_5=cv_6\cdot v_5.$  We get $c=0.$  Similarly, from $0=av_3\cdot v_2+bv_4\cdot v_2,$ we get $a=0.$  It turns out that $b=0.$ This is a contradiction.

On the other hand, let
\[A=\left[ \begin{matrix} 1 &1 &0&0&-1&0&-4\\
1&2&1&0&-2&0&-3\\
0&1&2&1&0&1&-1\\
0&0&1&2&0&3&-1\\
-1&-2&0&0&4&-1&0\\
0&0&1&3&-1&5&0\\
-4&-3&-1&-1&0&0&22\\
\end{matrix}
\right]\]

Let \[B=\left[\begin{matrix} 1 &1 &0&0&-1&0&-4\\
0&0&1&1&1&1&-2\\
0&0&0&1&-1&2&1\\
0&1&1&0&-1&0&1\\
\end{matrix}
\right]\]
Then $B^{T}B=A$, $A\in \P(G998)$ and $\rank(B)=4.$
\end{example}

\begin{example}
\[msr(G999)=4\]
The tree size of $G999$ is $5$.

Let

\[A=\left[ \begin{matrix} 1 &1 &0&0&0&1&1\\
1&2&1&0&0&1&0\\
0&1&5&2&0&2&-1\\
0&0&2&2&1&0&1\\
0&0&0&1&1&-1&1\\
1&1&2&0&-1&3&0\\
1&0&-1&1&1&0&3\\
\end{matrix}
\right]\]

Let \[B=\left[\begin{matrix} 1 &1 &0&0&0&1&1\\
0&0&2&1&0&1&0\\
0&0&0&1&1&-1&1\\
0&1&1&0&0&0&-1\\
\end{matrix}
\right]\]
Then $B^{T}B=A$, $A\in \P(G999)$ and $\rank(B)=4.$
\end{example}

\begin{example}
\[msr(G1000)=3\]
The tree size of $G1000$ is $4$.

Let \[B=\left[\begin{matrix} 1&-5&0&0&0&1&2\\
0&-2&1&1&0&-2&1\\
0&1&0&2&1&1&0\\
\end{matrix}
\right]\]
Let \[A=B^{T}B\]
Then $A\in \P(G1000)$ and $\rank(A)=3$.
\end{example}

\begin{example}
\[msr(G1003)=4\]
The tree size of $G1003$ is $5$.

Let \[B=\left[\begin{matrix} 1&1&0&0&0&1&2\\
0&0&2&1&0&1&-2\\
0&0&3&0&1&1&1\\
0&1&1&0&0&-1&1\\
\end{matrix}
\right]\]
Let \[A=B^{T}B\]
Then $A\in \P(G1003)$ and $\rank(A)=4$.

\end{example}

\begin{example}
\[msr(G1004)=4\]
The tree size of $G1004$ is $5$. 

Let \[B=\left[\begin{matrix} 1&1&0&0&0&1&3\\
0&0&1&1&0&1&-1\\
0&0&2&0&1&-1&1\\
0&1&1&0&0&1&-1\\
\end{matrix}
\right]\]
Let \[A=B^{T}B\]
Then $A\in \P(G1004)$  and $\rank(A)= 4$.

\end{example}

\begin{example}
\[msr(G1005)=3\]
The tree size of $G1005$ is $4$.

Let \[B=\left[\begin{matrix} 1&1&0&2&0&1&0\\
0&1&1&-2&1&0&0\\
0&0&0&1&2&-2&1\\
\end{matrix}
\right]\]
Let \[A=B^{T}B\]
Then $A\in \P(G1005)$ and $\rank(A)=3$.
\end{example}

\begin{example}
\[msr(G1006)=4\]
The tree size of $G1006$ is $4$. 
If $v=\{v_1,v_2,v_3, v_4, v_5,v_6, v_7\}$ is a representation of $G1006,$ then $v_1\perp\{v_3, v_4\},$ $v_2\perp\{v_4, v_7\},$ $v_3\perp\{v_5, v_6\},$ $v_4\perp v_5,$ $v_5\perp v_7,$ $v_6\perp v_7.$  We know $v_5\perp\{v_3,v_4, v_7\}.$ If $msr(G1006)=3,$ then there exists $a,b,c$ such that $av_3+bv_4+cv_7=0.$  From $av_3\cdot v_1+bv_4\cdot v_1+cv_7\cdot v_1=0,$ we get $c=0.$  From $av_3\cdot v_6+bv_4\cdot v_6=0,$ we get $b=0.$  Hence $av_3=0$  and it follows that $a=0.$  This is a contradiction.

Let \[B=\left[\begin{matrix} 1&1&0&0&1&1&0\\
0&1&1&1&0&0&1\\
0&-1&0&1&-1&0&1\\
0&0&0&0&1&2&1\\
\end{matrix}
\right]\]
Let \[A=B^{T}B\]
Then $A\in \P(G1006)$ and $\rank(A)= 4$.

\end{example}

\begin{example}
\[msr(G1007)=4\]
The tree size of $G1007$ is $5$.

Let \[B=\left[\begin{matrix} 1&1&0&0&0&1&1\\
0&2&1&0&0&-6&0\\
0&3&0&1&0&1&1\\
0&4&0&0&1&2&-1\\
\end{matrix}
\right]\]
Let \[A=B^{T}B\]
Then $A\in \P(G1007)$ and $\rank(A)= 4$.

\end{example}

\begin{example}
\[msr(G1053)=4\]
The tree size of $G1053$ is $5$.

Let \[B=\left[\begin{matrix} 1&1&0&0&0&1&-3\\
0&-2&1&1&0&2&1\\
0&1&1&2&0&1&1\\
0&1&0&0&1&1&0\\
\end{matrix}
\right]\]
Let \[A=B^{T}B\]
Then $A\in \P(G1053)$ and  $\rank(A)= 4$.

\end{example}

\begin{example}
\[msr(G1056)=3\]
The tree size of $G1056$ is $4$.

Let \[B=\left[\begin{matrix} 1&1&0&0&1&0&1\\
0&1&1&0&0&1&0\\
0&1&0&1&1&2&-1\\
\end{matrix}
\right]\]
Let \[A=B^{T}B\]
Then $A\in \P(G1056)$ and   $\rank(A)= 3$.
\end{example}

\begin{example}
\[msr(G1060)=3\]
The tree size of $G1060$ is $4$.

Let \[B=\left[\begin{matrix} 1&1&0&0&0&1&1\\
0&0&1&1&0&1&-1\\
0&1&1&0&1&0&0\\
\end{matrix}
\right]\]
Let \[A=B^{T}B\]
Then $A\in \P(G1060)$ and $\rank(A)= 3$.
\end{example}

\begin{example}
\[msr(G1065)=4\]
The tree size of $G1065$ is $4$.  If $v=\{v_1,v_2,v_3, v_4, v_5,v_6, v_7\}$ is a representation of $G1065,$ then $v_1\perp\{v_6\},$ $v_2\perp\{v_3,v_4, v_5,v_6\},$ $v_3\perp\{v_4, v_7\},$ $v_5\perp v_7.$    If $msr(G1065)=3,$ then there exists $a,b,c$ such that $v_5=av_2+bv_3+cv_4.$  From $av_2\cdot v_2+bv_3\cdot v_2+cv_4\cdot v_2=v_5\cdot v_2=0,$ we get $a=0.$  From $0=v_5\cdot v_7=bv_3\cdot v_7+cv_4\cdot v_7,$ we get $c=0.$  Hence $v_5=bv_3.$ It follows that $v_5\cdot v_4=bv_3\cdot v_4=0.$    This is a contradiction.

Let \[B=\left[\begin{matrix} 1&1&0&0&0&0&2\\
1&0&1&0&1&1&0\\
-2&0&0&1&1&1&1\\
1&0&0&0&-1&1&1\\
\end{matrix}
\right]\]
Let \[A=B^{T}B\]
Then $A\in \P(G1065)$ and $\rank(A)=4$.

\end{example}

\begin{example}
\[msr(G1069)=4\]
The tree size of $G1069$ is $4$.  If $v=\{v_1,v_2,v_3, v_4, v_5,v_6, v_7\}$ is a representation of $G1069,$ then $v_1\perp\{v_3,v_4\},$ $v_2\perp\{v_4, \},$ $v_3\perp\{v_7\}$,$v_4\perp v_7$, $v_5\perp \{v_6,v_7\}$, and $v_6\perp v_7.$    If $msr(G1069)=3,$ then  $v_1=kv_7.$  It follows that $v_1\cdot v_6=kv_7\cdot v_6=0.$  This is a contradiction.

Let \[B=\left[\begin{matrix} 1&1&0&0&0&0&1\\
0&-3&1&1&1&1&0\\
0&1&2&3&1&-2&0\\
1&0&0&0&1&1&0\\
\end{matrix}
\right]\]
Let \[A=B^{T}B\]
Then $A\in \P(G1069)$ and $\rank(A)=4$.

\end{example}

\begin{example}
\[msr(G1084)=4\]
The tree size of $G1084$ is $4$.  If $v=\{v_1,v_2,v_3, v_4, v_5,v_6, v_7\}$ is a representation of $G1084,$ then $v_1\perp\{v_3,v_4,v_5\},$ $v_2\perp\{v_5, \},$ $v_3\perp\{v_4,v_7\}$,$v_4\perp v_7$,  and $v_6\perp v_7.$    If $msr(G1084)=3,$ then  $v_1=kv_7.$  It follows that $v_1\cdot v_6=kv_7\cdot v_6=0.$  This is a contradiction.

Let \[B=\left[\begin{matrix} 0&1&1&0&1&1&0\\
0&1&0&1&-2&1&0\\
0&1&0&0&1&-1&1\\
1&1&0&0&0&1&1\\
\end{matrix}
\right]\]
Let \[A=B^{T}B\]
Then $A\in \P(G1084)$ and $rank(A)=4$.

\end{example}

\begin{example}
\[msr(G1089)=4\]
The tree size of $G1089$ is $5$.

Let \[B=\left[\begin{matrix} 1&1&0&0&0&-4&1\\
0&1&1&1&0&1&1\\
0&0&0&1&1&-1&-2\\
0&-1&0&1&0&1&1\\
\end{matrix}
\right]\]
Let \[A=B^{T}B\]
Then $A\in \P(G1089)$ and $\rank(A)=4$.
\end{example}

\begin{example}
\[msr(G1090)=3\]
The tree size of $G1090$ is $4$.

Let \[B=\left[\begin{matrix} 1&1&0&0&0&-1&2\\
0&-2&1&1&0&0&-2\\
0&1&0&2&1&2&1\\
\end{matrix}
\right]\]
Let \[A=B^{T}B\]
Then $A\in \P(G1090)$ and $\rank(A)= 3$.
\end{example}

\begin{example}
\[msr(G1091)=4\]
The tree size of $G1091$ is $4$.  If $v=\{v_1,v_2,v_3, v_4, v_5,v_6, v_7\}$ is a representation of $G1084,$ then $v_1\perp\{v_3,v_4,v_5\},$ $v_2\perp\{v_4,v_5, \},$ $v_3\perp\{v_7\}$,$v_5\perp v_7$,  and $v_6\perp v_7.$    If $msr(G1091)=3,$ then  $v_1=kv_7.$  It follows that $v_1\cdot v_6=kv_7\cdot v_6=0.$  This is a contradiction.

Let \[B=\left[\begin{matrix} 1&1&0&0&0&1&2\\
0&1&1&1&1&1&0\\
0&-1&0&1&1&1&1\\
0&0&0&0&1&3&-1\\
\end{matrix}
\right]\]
Let \[A=B^{T}B\]
Then $A\in \P(G1091)$ and $\rank(A)= 4$.

\end{example}

\begin{example}
\[msr(G1092)=4\]
The tree size of $G1092$ is $4$.  If $v=\{v_1,v_2,v_3, v_4, v_5,v_6, v_7\}$ is a representation of $G1084,$ then $v_1\perp\{v_3,v_4\},$ $v_2\perp\{v_4 \},$ $v_3\perp\{v_6,v_7\}$,$v_5\perp\{v_6, v_7\}$,  and $v_6\perp v_7.$    If $msr(G1092)=3,$ then  $v_5=kv_3.$  It follows that $v_5\cdot v_1=kv_3\cdot v_1=0.$  This is a contradiction.

Let \[B=\left[\begin{matrix} 0&1&1&1&1&0&0\\
1&1&0&1&0&1&0\\
-2&1&0&1&0&0&1\\
1&-3&0&1&1&0&0\\
\end{matrix}
\right]\]
Let \[A=B^{T}B\]
Then $A\in \P(G1092)$ and $\rank(A)= 4$.

\end{example}

\begin{example}
\[msr(G1093)=4\]
The tree size of $G1093$ is $5$.

Let \[B=\left[\begin{matrix} 1&2&0&0&0&1&1\\
0&1&1&1&0&2&0\\
0&0&0&1&1&1&1\\
0&-1&0&1&0&-2&1\\
\end{matrix}
\right]\]
Let \[A=B^{T}B\]
Then $A\in \P(G1093)$ and $\rank(A)= 4$.
\end{example}

\begin{example}
\[msr(G1094)=3\]
The tree size of $G1094$ is $4$.

Let \[B=\left[\begin{matrix} 1&1&0&0&0&2&-1\\
0&0&1&1&0&1&1\\
0&1&1&0&1&-1&-1\\
\end{matrix}
\right]\]
Let \[A=B^{T}B\]
Then $A\in \P(G1094)$ and $\rank(A)=3$.
\end{example}

\begin{example}
\[msr(G1095)=3\]
The tree size of $G1095$ is $4$.

Let \[B=\left[\begin{matrix} 1&1&1&0&0&1&0\\
1&0&0&1&1&1&0\\
-1&0&1&1&0&0&1\\
\end{matrix}
\right]\]
Let \[A=B^{T}B\]
Then $A\in \P(G1095)$ and $\rank(A)= 3$.
\end{example}

\begin{example}
\[msr(G1096)=3\]
The tree size of $G1096$ is $4$.

Let \[B=\left[\begin{matrix} 1&1&0&0&-3&0&1\\
0&0&1&1&2&0&1\\
0&1&1&0&-1&1&-1\\
\end{matrix}
\right]\]
Let \[A=B^{T}B\]
Then $A\in \P(G1096)$ and $\rank(A)= 3$.
\end{example}

\begin{example}
\[msr(G1097)=4\]
The tree size of $G1097$ is $4$. 
If $v=\{v_1,v_2,v_3, v_4, v_5,v_6, v_7\}$ is a representation of $G1097,$ then $v_1\perp\{v_3,v_4\},$ $v_2\perp\{v_4 ,v_5, v_6\},$ $v_3\perp\{v_5,v_6\}$, and $v_5\perp\{v_7\}$.    If $msr(G1097)=3,$ then  $v_2=kv_3.$  It follows that $v_2\cdot v_1=kv_3\cdot v_1=0.$  This is a contradiction.

Let \[B=\left[\begin{matrix} 1&1&1&0&0&0&1\\
0&0&1&1&1&1&-1\\
1&0&-1&0&1&1&1\\
0&0&0&0&0&1&1\\
\end{matrix}
\right]\]
Let \[A=B^{T}B\]
Then $A\in \P(G1097)$ and $\rank(A)= 4$.

\end{example}

\begin{example}
\[msr(G1100)=3\]
The tree size of $G1100$ is $4$.

Let \[B=\left[\begin{matrix} 1&1&0&0&1&-2&0\\
0&0&1&1&1&3&0\\
0&2&1&0&-1&1&1\\
\end{matrix}
\right]\]
Let \[A=B^{T}B\]
Then $A\in \P(G1100)$ and $\rank(A)= 3$.
\end{example}

\begin{example}
\[msr(G1101)=4\]
The tree size of $G1101$ is $4$.  If $v=\{v_1,v_2,v_3, v_4, v_5,v_6, v_7\}$ is a representation of $G1101,$ then $v_1\perp\{v_3,v_4,v_6\},$ $v_2\perp\{v_4, v_7\},$ $v_3\perp\{v_7\}$, $v_4\perp v_5,$ and $v_5\perp\{v_6\}$.    If $msr(G1101)=3,$ then  $v_5=kv_1.$  It follows that $v_5\cdot v_3=kv_3\cdot v_1=0.$  This is a contradiction.

Let \[B=\left[\begin{matrix} 1&1&0&0&1&0&1\\
0&1&1&1&1&1&0\\
0&-1&0&1&-1&1&2\\
0&1&0&0&0&1&1\\
\end{matrix}
\right]\]
Let \[A=B^{T}B\]
Then $A\in \P(G1101)$ and $\rank(A)= 4$.

\end{example}

\begin{example}
\[msr(G1104)=3\]
The tree size of $G1104$ is $4$.

Let \[B=\left[\begin{matrix} 1&1&0&0&1&1&-2\\
0&1&1&1&0&0&3\\
0&-1&0&1&2&1&1\\
\end{matrix}
\right]\]
Let \[A=B^{T}B\]
Then $A\in \P(G1104)$ and $\rank(A)= 3$.
\end{example}

\begin{example}
\[msr(G1105)=3\]
The tree size of $G1105$ is $4$.

Let \[B=\left[\begin{matrix} 5&1&1&1&0&0&0\\
-1&0&-3&1&1&1&0\\
4&0&-2&-1&0&1&1\\
\end{matrix}
\right]\]
Let \[A=B^{T}B\]
Then $A\in \P(G1105)$ and $\rank(A)= 3$.
\end{example}

\begin{example}
\[msr(G1123)=3\]
The tree size of $G1123$ is $4$.

Let \[B=\left[\begin{matrix} 0&1&1&0&0&1&0\\
1&2&0&1&1&0&0\\
-1&1&0&1&0&1&1\\
\end{matrix}
\right]\]
Let \[A=B^{T}B\]
Then $A\in \P(G1123)$ and $\rank(A)= 3$.
\end{example}

\begin{example}
\[msr(G1125)=3\]
The tree size of $G1125$ is $4$.

Let \[B=\left[\begin{matrix} 1&-2&0&1&1&0&1\\
0&1&1&1&0&0&0\\
0&1&0&0&1&1&2\\
\end{matrix}
\right]\]
Let \[A=B^{T}B\]
Then $A\in \P(G1125)$ and $\rank(A)= 3$.
\end{example}

\begin{example}
\[msr(G1135)=3\]
The tree size of $G1135$ is $4$.

Let \[B=\left[\begin{matrix} 1&1&0&0&0&0&1\\
0&-2&1&0&1&1&1\\
0&1&0&1&1&2&-1\\
\end{matrix}
\right]\]
Let \[A=B^{T}B\]
Then $A\in \P(G1135)$ and $\rank(A)= 3$.
\end{example}

\begin{example}
\[msr(G1145)=4\]
The tree size of $G1145$ is $4$. 
If $v=\{v_1,v_2,v_3, v_4, v_5,v_6, v_7\}$ is a representation of $G1145,$ then $v_1\perp\{v_3,v_4,v_5\},$ $v_2\perp\{v_4, v_5\},$ $v_4\perp\{v_6\}$, and $v_5\perp v_7$.    If $msr(G1145)=3,$ then  $v_2=kv_1.$  It follows that $v_2\cdot v_3=kv_3\cdot v_1=0.$  This is a contradiction.

Let \[B=\left[\begin{matrix} 1&1&0&0&0&1&1\\
0&0&1&1&1&0&1\\
0&0&1&0&1&1&-1\\
0&1&1&0&0&1&1\\
\end{matrix}
\right]\]
Let \[A=B^{T}B\]
Then $A\in \P(G1145)$ and $\rank(A)= 4$.

\end{example}

\begin{example}
\[msr(G1146)=3\]
The tree size of $G1146$ is $4$.

Let \[B=\left[\begin{matrix} 1&1&0&0&0&1&1\\
0&0&1&1&0&-1&1\\
0&1&1&0&1&2&0\\
\end{matrix}
\right]\]
Let \[A=B^{T}B\]
Then $A\in \P(G1146)$ and $\rank(A)= 3$.
\end{example}

\begin{example}
\[msr(G1148)=4\]
The tree size of $G1148$ is $5$. Let

\[A=\left[ \begin{matrix} 1 &0 &-1&0&0&1&1\\
0&1&-1&0&0&1&1\\
-1&-1&3&-1&0&-1&-3\\
0&0&-1&2&-1&-2&2\\
0&0&0&-1&1&1&-1\\
1&1&-1&-2&1&4&0\\
1&1&-3&2&-1&0&4\\
\end{matrix}
\right]\]
Let \[B=\left[\begin{matrix} 1&0&-1&0&0&1&1\\
0&1&-1&0&0&1&1\\
0&0&1&-1&0&1&-1\\
0&0&0&1&-1&-1&1\\
\end{matrix}
\right]\]
Then \[A=B^{T}B\]
 $A\in \P(G1148)$ and $\rank(A)= 4$.

\end{example}
\begin{example}
\[msr(G1149)=3\]
The tree size of $G1149$ is $4$. Let
 
\[A=\left[ \begin{matrix} 1 &-1 &0&0&1&1&0\\
-1&2&-1&0&-2&0&-1\\
0&-1&2&-1&2&0&1\\
0&0&-1&1&-1&-1&0\\
1&-2&2&-1&3&1&1\\
1&0&0&-1&1&3&-1\\
0&-1&1&0&1&-1&1\\
\end{matrix}
\right]\]
Let \[B=\left[\begin{matrix} 1&-1&0&0&1&1&0\\
0&1&-1&0&-1&1&-1\\
0&0&1&-1&1&1&0\\
\end{matrix}
\right]\]
Then \[A=B^{T}B\]
 $A\in \P(G1149)$ and $\rank(A)= 3$.

\end{example}
\begin{example}
\[msr(G1152)=3\]
The tree size of $G1152$ is $4$.  Let
\[A=\left[ \begin{matrix} 1 &1 &0&0&1&0&-2\\
1&3&-1&0&2&1&0\\
0&-1&1&1&2&0&-1\\
0&0&1&2&5&1&0\\
1&2&2&5&14&3&-1\\
0&1&0&1&3&1&1\\
-2&0&-1&0&-1&1&6\\
\end{matrix}
\right]\]
Let \[B=\left[\begin{matrix} 1&1&0&0&1&0&-2\\
0&1&-1&-1&-2&0&1\\
0&1&0&1&3&1&1\\
\end{matrix}
\right]\]
Then \[A=B^{T}B\]
 $A\in \P(G1152)$ and  $\rank(A)=3$.

\end{example}
\begin{example}
\[msr(G1153)=3\]
The tree size of $G1153$ is $4$.  Let
\[A=\left[ \begin{matrix} 1 &1 &0&0&-1&-1&1\\
1&2&-1&0&-4&-4&0\\
0&-1&5&2&7&7&-1\\
0&0&2&1&2&2&-1\\
-1&-4&7&2&14&7&0\\
-1&-4&7&2&7&14&0\\
1&0&-1&-1&0&0&3\\
\end{matrix}
\right]\]
Let \[B=\left[\begin{matrix} 1&1&0&0&-1&-1&1\\
0&1&-1&0&-3&-3&-1\\
0&0&2&1&2&2&-1\\
\end{matrix}
\right]\]
Then \[A=B^{T}B\]
 $A\in \P(G1153)$ and  $\rank(A)=3$.

\end{example}
\begin{example}
\[msr(G1154)=4\]
The tree size of $G1154$ is $4.$ First let's show $msr(G1154)>3.$

Let $\{v_1,v_2,v_3,v_4,v_5,v_6,v_7\}$ be a vector representation of $G1154.$
It follows that $v_1\perp v_3,$ $v_1\perp v_4,$ $v_1\perp v_7,$ $v_3\perp v_6,$ and $v_4\perp v_6.$ If $msr(G1154)=3,$ then $v_7\in \span \{v_3, v_4\}.$  It follows that $v_6\perp v_7.$  
We know that $v_6\cdot v_7\not= 0.$  This leads a contradiction.

On the other hand, let

\[A=\left[ \begin{matrix} 1 &1 &0&0&1&1&0\\
1&2&1&0&3&1&1\\
0&1&2&1&5&0&5\\
0&0&1&1&3&0&4\\
1&3&5&3&18&3&0\\
1&1&0&0&3&2&-7\\
0&1&5&4&0&-7&66\\
\end{matrix}
\right]\]
Let \[B=\left[\begin{matrix} 1&1&0&0&1&1&0\\
0&1&1&0&2&0&1\\
0&0&1&1&3&0&4\\
0&0&0&0&2&1&-7\\
\end{matrix}
\right]\]
Then \[A=B^{T}B\]
 $A\in \P(G1154)$ and  $\rank(A)=4$.
\end{example}
\begin{example}
\[msr(G1155)=3\]
The tree size of $G1155$ is $4$. 
Let
\[A=\left[ \begin{matrix} 1 &1 &0&0&-1&0&1\\
1&2&1&0&0&1&2\\
0&1&2&1&2&1&0\\
0&0&1&1&1&0&-1\\
-1&0&2&1&3&1&-1\\
0&1&1&0&1&1&1\\
1&2&0&-1&-1&1&3\\
\end{matrix}
\right]\]
Let \[B=\left[\begin{matrix} 1&1&0&0&-1&0&1\\
0&1&1&0&1&1&1\\
0&0&1&1&1&0&-1\\

\end{matrix}
\right]\]
Then \[A=B^{T}B\]
 $A\in \P(G1155)$ and  $\rank(A)=3$.

\end{example}
\begin{example}
\[msr(G1156)=3\]
The tree size of $G1156$ is $4$. 
Let
\[A=\left[ \begin{matrix} 1 &1 &0&0&1&1&0\\
1&2&1&0&0&3&-1\\
0&1&2&1&0&2&1\\
0&0&1&1&1&0&2\\
1&0&0&1&3&-1&3\\
1&3&2&0&-1&5&-2\\
0&-1&1&2&3&-2&5\\
\end{matrix}
\right]\]
Let \[B=\left[\begin{matrix} 1&1&0&0&1&1&0\\
0&1&1&0&-1&2&-1\\
0&0&1&1&1&0&2\\

\end{matrix}
\right]\]
Then \[A=B^{T}B\]
 $A\in \P(G1156)$ and  $\rank(A)=3$.

\end{example}
\begin{example}
\[msr(G1157)=3\]
The tree size of $G1157$ is $4$. 
Let
\[A=\left[ \begin{matrix} 1 &1 &0&0&1&0&-5\\
1&2&1&0&3&1&-3\\
0&1&2&1&3&-1&3\\
0&0&1&1&1&-2&1\\
1&3&3&1&6&0&0\\
0&1&-1&-2&0&5&0\\
-5&-3&3&1&0&0&30\\
\end{matrix}
\right]\]
Let \[B=\left[\begin{matrix} 1&1&0&0&1&0&-5\\
0&1&1&0&2&1&2\\
0&0&1&1&1&-2&1\\

\end{matrix}
\right]\]
Then \[A=B^{T}B\]
 $A\in \P(G1157)$ and  $\rank(A)=3$.

\end{example}
\begin{example}
\[msr(G1159)=3\]
The tree size of $G1159$ is $4$. Let

\[A=\left[ \begin{matrix} 1 &1 &0&0&2&1&2\\
1&5&2&0&0&-1&0\\
0&2&2&1&0&0&1\\
0&0&1&1&1&1&2\\
2&0&0&1&6&4&7\\
1&-1&0&1&4&3&5\\
2&0&1&2&7&5&9\\
\end{matrix}
\right]\]
Let \[B=\left[\begin{matrix} 1&1&0&0&2&1&2\\
0&2&1&0&-1&-1&-1\\
0&0&1&1&1&1&2\\

\end{matrix}
\right]\]
Then \[A=B^{T}B\]
 $A\in \P(G1159)$ and  $\rank(A)=3$.

\end{example}
\begin{example}
\[msr(G1160)=4\]
The tree size of $G1160$ is $5$. Let

\[A=\left[ \begin{matrix} 1 &1 &0&0&0&2&1\\
1&2&1&0&0&13&-1\\
0&1&3&3&0&7&-4\\
0&0&3&5&1&-2&-1\\
0&0&0&1&2&8&4\\
2&13&7&-2&8&165&0\\
1&-1&-4&-1&4&0&15\\
\end{matrix}
\right]\]
Let \[B=\left[\begin{matrix} 1&1&0&0&0&2&1\\
0&1&1&0&0&11&-2\\
0&0&1&2&1&2&1\\
0&0&1&1&-1&-6&-3\\
\end{matrix}
\right]\]
Then \[A=B^{T}B\]
 $A\in \P(G1160)$ and  $\rank(A)=4$.

\end{example}
\begin{example}
\[msr(G1165)=3\]
The tree size of $G1165$ is $4$. Let

\[A=\left[ \begin{matrix} 1 &2 &0&0&2&-1&1\\
2&6&-1&-1&3&-3&0\\
0&-1&1&0&0&-1&1\\
0&-1&0&1&1&2&1\\
2&3&0&1&5&0&3\\
-1&-3&-1&2&0&6&0\\
1&0&1&1&3&0&3\\
\end{matrix}
\right]\]
Let \[B=\left[\begin{matrix} 1&2&0&0&2&-1&1\\
0&1&-1&0&0&1&-1\\
0&1&0&-1&-1&-2&-1\\

\end{matrix}
\right]\]
Then \[A=B^{T}B\]
 $A\in \P(G1165)$ and  $\rank(A)=3$.

\end{example}
\begin{example}
\[msr(G1167)=3\]
The tree size of $G1167$ is $4$. Let

\[A=\left[ \begin{matrix} 1 &1 &0&0&3&1&1\\
1&6&-1&0&0&8&0\\
0&-1&2&3&0&-2&-1\\
0&0&3&5&-1&-1&-2\\
3&0&0&-1&11&-1&4\\
1&8&-2&-1&-1&11&0\\
1&0&-1&-2&4&0&2\\
\end{matrix}
\right]\]
Let \[B=\left[\begin{matrix} 1&1&0&0&3&1&1\\
0&2&-1&-1&-1&3&0\\
0&1&1&2&-1&1&-1\\

\end{matrix}
\right]\]
Then \[A=B^{T}B\]
 $A\in \P(G1167)$ and  $\rank(A)=3$.

\end{example}
\begin{example}
\[msr(G1168)=3\]
The tree size of $G1168$ is $4$. Let
 
\[A=\left[ \begin{matrix} 1 &1 &0&0&1&1&0\\
1&6&-1&0&4&0&-4\\
0&-1&2&3&0&-1&5\\
0&0&3&5&1&-2&7\\
1&4&0&1&3&0&-1\\
1&0&-1&-2&0&2&-2\\
0&-4&5&7&-1&-2&13\\
\end{matrix}
\right]\]
Let \[B=\left[\begin{matrix} 1&1&0&0&1&1&0\\
0&1&1&2&1&-1&2\\
0&2&-1&-1&1&0&-3\\

\end{matrix}
\right]\]
Then \[A=B^{T}B\]
 $A\in \P(G1168)$ and  $\rank(A)=3$.

\end{example}
\begin{example}
\[msr(G1170)=3\]
The tree size of $G1170$ is $4$. Let

\[A=\left[ \begin{matrix} 1 &1 &0&0&-1&-2&-2\\
1&2&1&0&0&-3&-1\\
0&1&2&1&2&0&-2\\
0&0&1&1&1&1&-3\\
-1&0&2&1&3&2&0\\
-2&-3&0&1&2&6&0\\
-2&-1&-2&-3&0&0&14\\
\end{matrix}
\right]\]
Let \[B=\left[\begin{matrix} 1&1&0&0&-1&-2&-2\\
0&1&1&0&1&-1&1\\
0&0&1&1&1&1&-3\\

\end{matrix}
\right]\]
Then \[A=B^{T}B\]
 $A\in \P(G1170)$ and  $\rank(A)=3$.

\end{example}

\begin{example}
\[msr(G1176)=3\]
The tree size of $G1176$ is $4$.

Let \[B=\left[\begin{matrix} 1&1&0&0&0&1&1\\
0&1&1&1&0&0&-1\\
0&1&0&1&1&1&-1\\
\end{matrix}
\right]\]
Let \[A=B^{T}B\]
$A\in \P(G1176)$ and  $\rank(A)=3$.

\end{example}
\begin{example}
\[msr(G1179)=3\]
The tree size of $G1179$ is $4$. Let

\[A=\left[ \begin{matrix} 1 &5 &0&0&1&3&6\\
5&45&2&0&7&5&0\\
0&2&2&3&2&-1&-3\\
0&0&3&5&3&0&0\\
1&7&2&3&3&2&3\\
3&5&-1&0&2&14&33\\
6&0&-3&0&3&33&81\\
\end{matrix}
\right]\]

Let \[B=\left[\begin{matrix} 1&5&0&0&1&3&6\\
0&2&-1&-2&-1&-1&-3\\
0&4&1&1&1&-2&-6\\
\end{matrix}
\right]\]
Then \[A=B^{T}B\]
$A\in \P(G1179)$ and  $\rank(A)=3$.

\end{example}

\begin{example}
\[msr(G1189)=3\]
The tree size of $G1189$ is $4$. Let

\[A=\left[ \begin{matrix} 1 &1 &0&0&1&0&1\\
1&6&3&1&2&7&5\\
0&3&2&0&0&5&3\\
0&1&0&2&2&-1&-1\\
1&2&0&2&3&-1&0\\
0&7&5&-1&-1&13&8\\
1&5&3&-1&0&8&6\\
\end{matrix}
\right]\]
Let \[B=\left[\begin{matrix} 1&1&0&0&1&0&1\\
0&2&1&1&1&2&1\\
0&1&1&-1&-1&3&2\\
\end{matrix}
\right]\]
Then \[A=B^{T}B\]
$A\in \P(G1189)$ and  $\rank(A)=3$.

\end{example}

\begin{example}
\[msr(G1191)=3\]
The tree size of $G1191$ is $4$. Let

\[A=\left[ \begin{matrix} 1 &1 &0&0&1&1&1\\
1&14&-1&5&3&8&2\\
0&-1&2&0&1&1&-2\\
0&5&0&2&1&3&0\\
1&3&1&1&2&3&0\\
1&8&1&3&3&6&0\\
1&2&-2&0&0&0&3\\
\end{matrix}
\right]\]
Let \[B=\left[\begin{matrix} 1&1&0&0&1&1&1\\
0&2&1&1&1&2&-1\\
0&3&-1&1&0&1&1\\
\end{matrix}
\right]\]
Then \[A=B^{T}B\]
$A\in \P(G1191)$ and  $\rank(A)=3$.

\end{example}
\begin{example}
\[msr(G1194)=3\]
The tree size of $G1194$ is $4$. Let

\[A=\left[ \begin{matrix} 1 &1 &0&0&1&0&1\\
1&6&-1&0&-2&-1&6\\
0&-1&2&3&0&5&-1\\
0&0&3&5&-1&8&0\\
1&-2&0&-1&3&-1&-2\\
0&-1&5&8&-1&13&-1\\
1&6&-1&0&-2&-1&6\\
\end{matrix}
\right]\]
Let \[B=\left[\begin{matrix} 1&1&0&0&1&0&1\\
0&2&-1&-1&-1&-2&2\\
0&1&1&2&-1&3&1\\
\end{matrix}
\right]\]
Then \[A=B^{T}B\]
$A\in \P(G1194)$ and  $\rank(A)=3$.

\end{example}

\begin{example}
\[msr(G1195)=3\]
The tree size of $G1195$ is $4$. Let

\[A=\left[ \begin{matrix} 1 &1 &0&0&4&0&2\\
1&3&1&0&5&2&6\\
0&1&5&3&11&1&-6\\
0&0&3&2&7&0&-1\\
4&5&11&7&41&1&0\\
0&2&1&0&1&2&-9\\
2&-7&-6&-1&0&-9&45\\
\end{matrix}
\right]\]
Let \[B=\left[\begin{matrix} 1&1&0&0&4&0&2\\
0&1&2&1&4&1&-5\\
0&1&-1&-1&-3&1&-4\\
\end{matrix}
\right]\]
Then \[A=B^{T}B\]
$A\in \P(G1195)$ and  $\rank(A)=3$.

\end{example}
\begin{example}
\[msr(G1196)=3\]
The tree size of $G1196$ is $4$. Let

\[A=\left[ \begin{matrix} 1 &1 &0&0&1&2&2\\
1&19&9&0&-2&-1&8\\
0&9&5&-1&0&0&3\\
0&0&-1&2&-3&-3&0\\
1&-2&0&-3&5&7&1\\
2&-1&0&-3&7&9&3\\
2&8&3&0&1&3&6\\
\end{matrix}
\right]\]
Let \[B=\left[\begin{matrix} 1&1&0&0&1&2&2\\
0&3&1&1&-2&-2&1\\
0&3&2&-1&1&1&1\\
\end{matrix}
\right]\]
Then \[A=B^{T}B\]
$A\in \P(G1196)$ and  $\rank(A)=3$.

\end{example}
\begin{example}
\[msr(G1197)=3\]
The tree size of $G1197$ is $4$. Let

\[A=\left[ \begin{matrix} 1 &1 &0&0&1&1&2\\
1&2&-1&-2&-2&0&0\\
0&-1&5&0&7&-1&0\\
0&-2&0&5&4&3&5\\
1&-2&7&4&14&2&6\\
1&0&-1&3&2&3&5\\
2&0&0&5&6&5&9\\
\end{matrix}
\right]\]
Let \[B=\left[\begin{matrix} 1 &1 &0&0&1&1&2\\
0&-1&1&2&3&1&2\\
0&0&2&-1&2&-1&-1\\
\end{matrix}
\right]\]
Then $B^{T}B=A$, $A\in \P(G1197)$ and  $\rank(A)=3$.

\end{example}
\begin{example}
\[msr(G1199)=3\]
The tree size of $G1199$ is $4$. Let

\[A=\left[ \begin{matrix} 1 &1 &0&1&0&2&0\\
1&14&5&13&22&9&5\\
0&5&2&5&9&1&1\\
1&13&5&14&23&0&0\\
0&22&9&23&41&0&2\\
2&9&1&0&0&45&23\\
0&5&1&0&2&23&13\\
\end{matrix}
\right]\]
Let \[B=\left[\begin{matrix} 1 &1 &0&1&0&2&0\\
0&2&1&3&5&-4&-2\\
0&3&1&2&4&5&3\\
\end{matrix}
\right]\]
Then $B^{T}B=A$, $A\in \P(G1199)$ and  $\rank(A)=3$.

\end{example}

\begin{example}
\[msr(G1200)=3\]
The tree size of $G1200$ is $4$. Let

\[A=\left[ \begin{matrix} 1 &1 &0&0&1&0&1\\
1&3&-1&0&6&1&1\\
0&-1&5&3&-4&7&3\\
0&0&3&2&-1&5&2\\
1&6&-4&-1&14&0&0\\
0&1&7&5&0&13&5\\
1&1&3&2&0&5&3\\
\end{matrix}
\right]\]
Let \[B=\left[\begin{matrix} 1 &1 &0&0&1&0&1\\
0&-1&2&1&-3&2&1\\
0&1&1&1&2&3&1\\
\end{matrix}
\right]\]
Then $B^{T}B=A$, $A\in \P(G1200)$ and  $\rank(A)=3$.

\end{example}
\begin{example}
\[msr(G1202)=3\]
The tree size of $G1202$ is $4$. Let

\[A=\left[ \begin{matrix} 1 &1 &0&0&0&1&2\\
1&6&-1&3&0&9&-1\\
0&-1&2&0&3&-1&9\\
0&3&0&2&1&5&1\\
0&0&3&1&5&1&14\\
1&9&-1&5&1&14&0\\
2&-1&9&1&14&0&45\\
\end{matrix}
\right]\]
Let \[B=\left[\begin{matrix} 1 &1 &0&0&0&1&2\\
0&-2&1&-1&1&-3&4\\
0&1&1&1&2&2&5\\
\end{matrix}
\right]\]
Then $B^{T}B=A$, $A\in \P(G1202)$ and  $\rank(A)=3$.

\end{example}
\begin{example}
\[msr(G1205)=3\]
The tree size of $G1205$ is $4$. Let

\[A=\left[ \begin{matrix} 1 &1 &0&0&1&-1&2\\
1&6&3&0&-1&0&7\\
0&3&5&4&2&3&3\\
0&0&4&5&4&3&0\\
1&-1&2&4&5&1&0\\
-1&0&3&3&1&3&-1\\
2&7&3&0&0&-1&9\\
\end{matrix}
\right]\]
Let \[B=\left[\begin{matrix} 1 &1 &0&0&1&-1&2\\
0&-1&1&2&2&1&-1\\
0&2&2&1&0&1&2\\
\end{matrix}
\right]\]
Then $B^{T}B=A$, $A\in \P(G1205)$ and  $\rank(A)=3$.

\end{example}

\begin{example}
\[msr(G1207)=3\]
The tree size of $G1207$ is $4$. Let

\[A=\left[ \begin{matrix} 1 &1 &0&0&-1&1&1\\
1&6&3&0&0&5&-3\\
0&3&5&4&3&8&0\\
0&0&4&5&3&7&3\\
-1&0&3&3&3&4&0\\
1&5&8&7&4&14&2\\
1&-3&0&3&0&2&6\\
\end{matrix}
\right]\]
Let \[B=\left[\begin{matrix} 1 &1 &0&0&-1&1&1\\
0&-1&1&2&1&2&2\\
0&2&2&1&1&3&-1\\
\end{matrix}
\right]\]
Then $B^{T}B=A$, $A\in \P(G1207)$ and  $\rank(A)=3$.
 
\end{example}
\begin{example}
\[msr(G1208)=3\]
The tree size of $G1208$ is $4$. Let

\[A=\left[ \begin{matrix} 1 &1 &0&0&1&1&1\\
1&3&3&1&0&0&3\\
0&3&5&0&-1&-3&4\\
0&1&0&5&-2&4&-2\\
1&0&-1&-2&2&0&1\\
1&0&-3&4&0&6&-3\\
1&3&4&-2&1&-3&5\\
\end{matrix}
\right]\]
Let \[B=\left[\begin{matrix} 1 &1 &0&0&1&1&1\\
0&1&2&-1&0&-2&2\\
0&1&1&2&-1&1&0\\
\end{matrix}
\right]\]
Then $B^{T}B=A$, $A\in \P(G1208)$ and  $\rank(A)=3$.

\end{example}

\begin{example}
\[msr(G1209)=3\]
The tree size of $G1209$ is $4$. Let

\[A=\left[ \begin{matrix} 1 &1 &0&0&2&1&1\\
1&6&3&0&1&7&5\\
0&3&5&4&1&2&0\\
0&0&4&5&2&-2&-3\\
2&1&1&2&5&0&0\\
1&7&2&-2&0&9&7\\
1&5&0&-3&0&7&6\\
\end{matrix}
\right]\]
Let \[B=\left[\begin{matrix} 1 &1 &0&0&2&1&1\\
0&-1&1&2&1&-2&-2\\
0&2&2&1&0&2&1\\
\end{matrix}
\right]\]
Then $B^{T}B=A$, $A\in \P(G1209)$ and  $\rank(A)=3$.

\end{example}

\begin{example}
\[msr(G1210)=3\]
The tree size of $G1210$ is $4$. Let

\[A=\left[ \begin{matrix} 1 &1 &0&0&1&3&1\\
1&3&-1&0&2&2&9\\
0&-1&5&3&7&-1&-7\\
0&0&3&2&5&-1&-2\\
1&2&7&5&14&0&-1\\
3&2&-1&-1&0&10&0\\
1&9&-7&-2&-1&0&35\\
\end{matrix}
\right]\]
Let \[B=\left[\begin{matrix} 1 &1 &0&0&1&3&1\\
0&-1&2&1&2&0&-5\\
0&1&1&1&3&-1&3\\
\end{matrix}
\right]\]
Then $B^{T}B=A$, $A\in \P(G1210)$ and  $\rank(A)=3$.

\end{example}

\begin{example}
\[msr(G1211)=2\]
The tree size of $G1211$ is $3$. Let

\[A=\left[ \begin{matrix} 1 &1 &2&-1&-1&-2&0\\
1&5&0&-1&-5&0&2\\
2&0&5&-2&0&-5&-1\\
-1&-1&-2&1&1&2&0\\
-1&-5&0&1&5&0&-2\\
-2&0&-5&2&0&5&1\\
0&2&-1&0&-2&1&1\\
\end{matrix}
\right]\]
Let \[B=\left[\begin{matrix} 1 &1 &2&-1&-1&-2&0\\
0&2&-1&0&-2&1&1\\
\end{matrix}
\right]\]
Then $B^{T}B=A$, $A\in \P(G1211)$ and  $\rank(A)=2$.

\end{example}

\begin{example}
\[msr(G1212)=3\]
The tree size of $G1212$ is $4$. Let

\[A=\left[ \begin{matrix} 1 &1 &0&0&3&1&2\\
1&3&1&0&0&2&3\\
0&1&5&3&0&8&-1\\
0&0&3&2&1&5&-1\\
3&0&0&1&10&4&4\\
1&2&8&5&4&14&0\\
2&3&-1&-1&4&0&5\\
\end{matrix}
\right]\]
Let \[B=\left[\begin{matrix} 1 &1 &0&0&3&1&2\\
0&-1&1&1&2&2&-1\\
0&1&2&1&-1&3&0\\
\end{matrix}
\right]\]
Then $B^{T}B=A$, $A\in \P(G1212)$ and  $\rank(A)=3$.

\end{example}
\begin{example}
\[msr(G1222)=3\]
The tree size of $G1222$ is $3.$ First let's show $msr(G1222)>2.$

If $v=\{v_1,v_2,v_3, v_4, v_5,v_6, v_7\}$ is a representation of $G1222,$ 
It follows that $v_1\perp v_3,$ $v_3\perp v_4,$ $v_4\perp v_6.$  If $msr(G1222)=2,$ then $v_1\perp  v_6.$
We know that $v_1\cdot v_6\not= 0.$  This leads a contradiction.

On the other hand, let

\[A=\left[ \begin{matrix} 1 &1 &0&1&1&1&1\\
1&3&2&1&6&2&1\\
0&2&2&0&5&1&0\\
1&1&0&3&2&0&-1\\
1&6&5&2&14&3&0\\
1&2&1&0&3&2&2\\
1&1&0&-1&0&2&3\\
\end{matrix}
\right]\]
Let \[B=\left[\begin{matrix} 1 &1 &0&1&1&1&1\\
0&1&1&-1&2&1&1\\
0&1&1&1&3&0&-1\\
\end{matrix}
\right]\]
Then $B^{T}B=A$, $A\in \P(G1222)$ and  $\rank(A)=3$.

\end{example}

\begin{example}
\[msr(G1224)=3\]
The tree size of $G1224$ is $3.$ First let's show $msr(G1224)>2.$

If $v=\{v_1,v_2,v_3, v_4, v_5,v_6, v_7\}$ is a representation of $G1224,$ 
It follows that $v_1\perp v_5,$ $v_1\perp v_6,$  If $msr(G1224)=2,$ then $v_5\cdot  v_6\not=0$
We know that $v_5\cdot v_6= 0.$  This leads a contradiction.

On the other hand, let

\[A=\left[ \begin{matrix} 1 &1 &0&1&0&0&1\\
1&5&2&3&2&-4&-1\\
0&2&2&2&3&-1&-1\\
1&3&2&3&3&-1&0\\
0&2&3&3&5&0&-1\\
0&-4&-1&-1&0&5&2\\
1&-1&-1&0&-1&2&2\\
\end{matrix}
\right]\]
Let \[B=\left[\begin{matrix} 1 &1 &0&1&0&0&1\\
0&2&1&1&1&-2&-1\\
0&0&1&1&2&1&0\\
\end{matrix}
\right]\]
Then $B^{T}B=A$, $A\in \P(G1224)$ and  $\rank(A)=3$.

\end{example}

\begin{example}
\[msr(G1228)=3\]
The tree size of $G1228$ is $3.$ First let's show $msr(G1228)>2.$

Let $v=\{v_1,v_2,v_3, v_4, v_5,v_6, v_7\}$ be a representation of $G1228.$ 
It follows that $v_1\perp v_3,$ $v_3\perp v_4,$ $v_4\perp v_7.$  If $msr(G1228)=2,$ then $v_1\perp v_7,$
We know that $v_1\cdot v_7\not= 0.$  This leads a contradiction.

On the other hand, let

\[A=\left[ \begin{matrix} 1 &1 &0&1&0&1&1\\
1&5&2&-1&2&5&3\\
0&2&2&2&3&-1&-1\\
1&3&2&0&3&1&1\\
1&-1&0&3&1&-2&0\\
1&5&1&-2&0&6&3\\
1&3&1&0&1&3&2\\
\end{matrix}
\right]\]
Let \[B=\left[\begin{matrix} 1 &1 &0&1&0&1&1\\
0&2&1&-1&1&2&1\\
0&0&1&1&2&-1&0\\
\end{matrix}
\right]\]
Then $B^{T}B=A$, $A\in \P(G1228)$ and  $\rank(A)=3$.
 
\end{example}

\begin{example}
\[msr(G1230)=2\]
The tree size of $G1230$ is $3$. Let
 
\[A=\left[ \begin{matrix} 1 &1 &0&1&-1&1&-1\\
1&5&2&3&-3&-1&1\\
0&2&1&1&-1&-1&1\\
1&3&1&2&-2&0&0\\
-1&-3&-1&-2&2&0&0\\
1&-1&-1&0&0&2&-2\\
-1&1&1&0&0&-2&2\\
\end{matrix}
\right]\]
Let \[B=\left[\begin{matrix} 1 &1 &0&1&-1&1&-1\\
0&2&1&1&-1&-1&1\\

\end{matrix}
\right]\]
Then $B^{T}B=A$, $A\in \P(G1230)$ and  $\rank(A)=2$.
 
\end{example}

\begin{example}
\[msr(G1231)=3\]
The tree size of $G1231$ is $3.$ First let's show $msr(G1231)>2.$

Let $v=\{v_1,v_2,v_3, v_4, v_5,v_6, v_7\}$ be a representation of $G1231.$ It follows that $v_2\perp v_4,$ $v_2\perp v_7,$ $v_5\perp v_7.$  If $msr(G1231)=2,$ then $v_4\perp v_5.$
We know that $v_4\cdot v_5\not= 0.$  This leads a contradiction.

On the other hand, let

\[A=\left[ \begin{matrix} 1 &1 &0&1&1&1&1\\
1&2&1&0&4&2&0\\
0&1&5&-1&7&3&1\\
1&0&-1&2&-2&0&2\\
1&4&7&-2&14&6&0\\
1&2&3&0&6&3&1\\
1&0&1&2&0&1&3\\
\end{matrix}
\right]\]
Let \[B=\left[\begin{matrix} 1 &1 &0&1&1&1&1\\
0&1&1&-1&3&1&-1\\
0&0&2&0&2&1&1\\
\end{matrix}
\right]\]
Then $B^{T}B=A$, $A\in \P(G1231)$ and  $\rank(A)=3$.

\end{example}

\begin{example}
\[msr(G1233)=3\]
The tree size of $G1233$ is $3.$ First let's show $msr(G1233)>2.$

Let $v=\{v_1,v_2,v_3, v_4, v_5,v_6, v_7\}$ be a representation of $G1233.$ It follows that $v_2\perp v_4,$ $v_2\perp v_6.$  If $msr(G1233)=2,$ then $v_4\cdot  v_6\not=0$
We know that $v_4\cdot v_6= 0.$  This leads a contradiction.

On the other hand, let

\[A=\left[ \begin{matrix} 1 &1 &0&1&1&-1&-3\\
1&2&1&0&1&0&-3\\
0&1&5&1&6&5&2\\
1&0&1&3&4&0&-2\\
1&1&6&4&10&5&0\\
-1&0&5&0&5&6&5\\
-3&-3&2&-2&0&5&10\\
\end{matrix}
\right]\]
Let \[B=\left[\begin{matrix} 1 &1 &0&1&1&-1&-3\\
0&1&1&-1&0&1&0\\
0&0&2&1&3&2&1\\
\end{matrix}
\right]\]
Then $B^{T}B=A$, $A\in \P(G1233)$ and  $\rank(A)=3$.
 
\end{example}

\begin{example}
\[msr(G1241)=3\]
The tree size of $G1241$ is $4$. Let

\[A=\left[ \begin{matrix} 1 &1 &0&0&1&1&-1\\
1&6&1&0&2&3&-2\\
0&1&2&3&2&1&1\\
0&0&3&5&3&1&2\\
1&2&2&3&3&2&0\\
1&3&1&1&2&2&-1\\
-1&-2&1&2&0&-1&2\\
\end{matrix}
\right]\]
Let \[B=\left[\begin{matrix} 1 &1 &0&0&1&1&-1\\
0&2&1&1&1&1&0\\
0&-1&1&2&1&0&1\\
\end{matrix}
\right]\]
Then $B^{T}B=A$, $A\in \P(G1241)$ and  $\rank(A)=3$.

\end{example}

\begin{example}
\[msr(G1242)=2\]
The tree size of $G1242$ is $3$. Let

\[A=\left[ \begin{matrix} 1 &1 &0&-1&1&0&1\\
1&2&1&-2&3&-1&0\\
0&1&1&-1&2&-1&-1\\
-1&-2&-1&2&-3&1&0\\
1&3&2&-3&5&-2&-1\\
0&-1&-1&1&-2&1&1\\
1&0&-1&0&-1&1&2\\
\end{matrix}
\right]\]
Let \[B=\left[\begin{matrix} 1 &1 &0&-1&1&0&1\\
0&1&1&-1&2&-1&-1\\
\end{matrix}
\right]\]
Then $B^{T}B=A$, $A\in \P(G1242)$ and  $\rank(A)=2$.

\end{example}

\begin{example}
\[msr(G1248)=2\]
The tree size of $G1248$ is $3$. Let

\[A=\left[ \begin{matrix} 1 &1 &0&1&1&1&2\\
1&10&3&4&7&-2&-1\\
0&3&1&1&2&-1&-1\\
1&4&1&2&3&0&1\\
1&7&2&3&5&-1&0\\
1&-2&-1&0&-1&2&3\\
2&-1&-1&1&0&3&5\\
\end{matrix}
\right]\]
Let \[B=\left[\begin{matrix} 1 &1 &0&1&1&1&2\\
0&3&1&1&2&-1&-1\\
\end{matrix}
\right]\]
Then $B^{T}B=A$, $A\in \P(G1248)$ and  $\rank(A)=2$.

\end{example}

\begin{example}
\[msr(G1250)=2\]
The tree size of $G1250$ is $3$. Let

\[A=\left[ \begin{matrix} 1 &1 &0&1&1&2&1\\
1&2&1&3&0&3&4\\
0&1&1&2&-1&1&3\\
1&3&2&5&-1&4&7\\
1&0&-1&-1&2&1&-2\\
2&3&1&4&1&5&5\\
1&4&3&7&-2&5&10\\
\end{matrix}
\right]\]
Let \[B=\left[\begin{matrix} 1 &1 &0&1&1&2&1\\
0&1&1&2&-1&1&3\\
\end{matrix}
\right]\]
Then $B^{T}B=A$, $A\in \P(G1250)$ and  $\rank(A)=2$.

\end{example}

\end{document}